\documentclass[10pt, twoside, twocolumn]{article}

\usepackage{amsmath,amsthm}

\usepackage[utf8]{inputenc}

\usepackage{enumerate}

\usepackage{fancyhdr}
\usepackage{cite}
\usepackage{enumitem}

\usepackage{graphics}
\usepackage{graphicx}
\usepackage{mathptmx}      
\usepackage{amssymb}
\usepackage{amsmath}
\usepackage[T1]{fontenc}
\usepackage{float}
\usepackage{subcaption}
\usepackage{xcolor}
\usepackage[utf8]{inputenc}
\inputencoding{latin1}
\inputencoding{utf8}
\usepackage{resizegather}
\usepackage{blindtext}
\usepackage{multicol}
\usepackage{array}
\newcolumntype{P}[1]{>{\centering\arraybackslash}p{#1}}
\usepackage{mwe}
\usepackage{hyperref}
\hypersetup{%
    colorlinks=True, 
    citecolor=blue,
    linkcolor=black}

  \theoremstyle{definition}

  \newtheorem*{assumption*}{Assumption}

  \numberwithin{equation}{section}
    
\usepackage{authblk}
\usepackage[T1]{fontenc}
\usepackage[utf8]{inputenc}

\author[1,2]{Thomas Gauthey\thanks{E-Mail:thomas.gauthey(at)geeps.centralesupelec.fr}}
\author[2]{Peter Gangl\thanks{E-Mail: gangl(at)math.tugraz.at}}
\author[1]{Maya Hage Hassan\thanks{E-Mail: maya.hage-hassan(at)centralesupelec.fr}}

\affil[1]{ Université Paris-Saclay, CentraleSupélec, CNRS, Laboratoire de Génie Electrique et Electronique de Paris, 91192, Gif-sur-Yvette, France. \\

\quad \quad Sorbonne Université, CNRS, Laboratoire de Génie Electrique et Electronique de Paris, 75252, Paris, France.}
\affil[2]{Technische Universität Graz, Institut für Angewandte Mathematik, 8010 Graz, Austria.}

\title{Multi-Material Topology Optimization with Continuous Magnetization Direction for Permanent Magnet Synchronous Reluctance Motors}

\date{\today}

\begin{document}

\maketitle

\begin{abstract}
Permanent magnet-assisted synchronous reluctance motors (PMSynRM) have a significantly higher average torque than synchronous reluctance motors. Thus, determining an optimal design results in a multi-material topology optimization problem, where one seeks to distribute ferromagnetic material, air and permanent magnets within the rotor in an optimal manner.This study proposed a novel density-based distribution scheme, which allows for continuous magnetization direction instead of a finite set of angles. Thus, an interpolation scheme is established between properties pertaining to magnets and non-linear materials. This allows for new designs to emerge without introducing complex geometric parameterization or relying on the user’s biases and intuitions. Toward reducing computation time, Nitsche-type mortaring is applied, allowing for free rotation of the rotor mesh relative to the stator mesh. The average torque is approximated using only four-point static positions. This study investigates several interpolation schemes and presents a new one inspired by the topological derivative. We propose to filter the final design for the magnetization angle using K-mean clustering accounting for technical feasibility constraints of magnets. Finally, the design of the electrical motor is proposed to maximize torque value.
\end{abstract}


\section{Introduction}
\label{intro}
Synchronous reluctance machines (SynRM) are standard in households and industrial applications, thanks to their cheap cost compared to permanent magnet motors and advances in manufacturing techniques. Although the deployment of these machines continues \cite{heidari_review_2021}, PMSynRM offers an excellent alternative for both structures, solving for SynRM, its poor power factor and, for permanent magnet machines (PMM), its cost. The design of these machines using parametric optimization often necessitates either complex analytical models or the use of Finite Element Analysis (FEA) relying heavily on experienced engineers and known good designs \cite{maroufian_pm_2017}. \\
Density based optimizations allow for bypassing such cumbersome frameworks. Although they were first developed for two materials application in continuous mechanics \cite{sigmund_99_2001}, a rise in n-materials optimization in the field of electromagnetics has allowed for new PMM and SynRM to emerge \cite{kim_magnetic_2010,guo_multimaterial_2020}.
In most optimizations where permanent magnets are involved, magnetization direction are fixed \cite{risticevic_design_2016} or limited to a set of a couple values \cite{lee_simultaneous_2011,choi_optimal_2010,choi_optimization_2010}. If continuous directions are considered during the optimization process, the final design is filtered to take into account only a couple of predefined directions to meet manufacturing constraints \cite{wang_topology_2005,ishikawa_topology_2015}. \\
In this paper we propose a simultaneous density-based optimization scheme consisting of three-material (air-iron-magnet) with a continuous magnetization direction. The proposition is applied to design the rotor of a distributed winding stator as described in \cite{mellak_synchronous_2018,gangl_multi-objective_2020} to maximize the mean torque under constraints. The final designs are filtered using an unbiased K-means heuristic for accounting for feasibility constraints.
Here, we propose also to accelerate the torque calculation through a four-point method.

\section{Problem description}
 We chose to investigate a SynRM described  in \cite{gangl_multi-objective_2020,mellak_synchronous_2018}, of which the rotor design had proven to be a challenging problem for topology optimization and use it for our PMSynRM optimization problem. 
\subsection{Geometry description}
\begin{figure}
\includegraphics[width=9 cm]{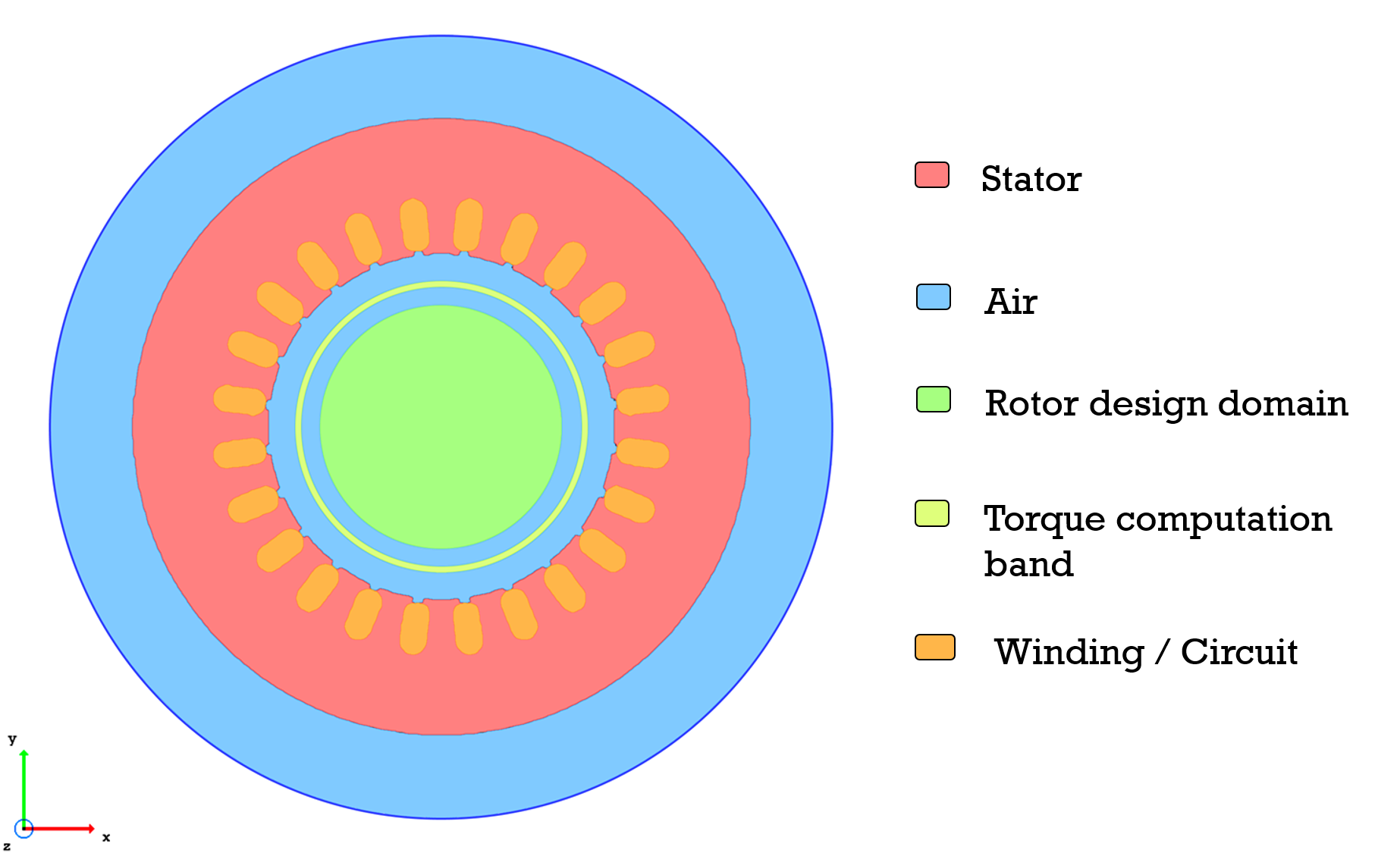}
\caption{Machine geometry }
\label{fig:machgeom}
\end{figure}

\begin{figure}
\includegraphics[width=9 cm]{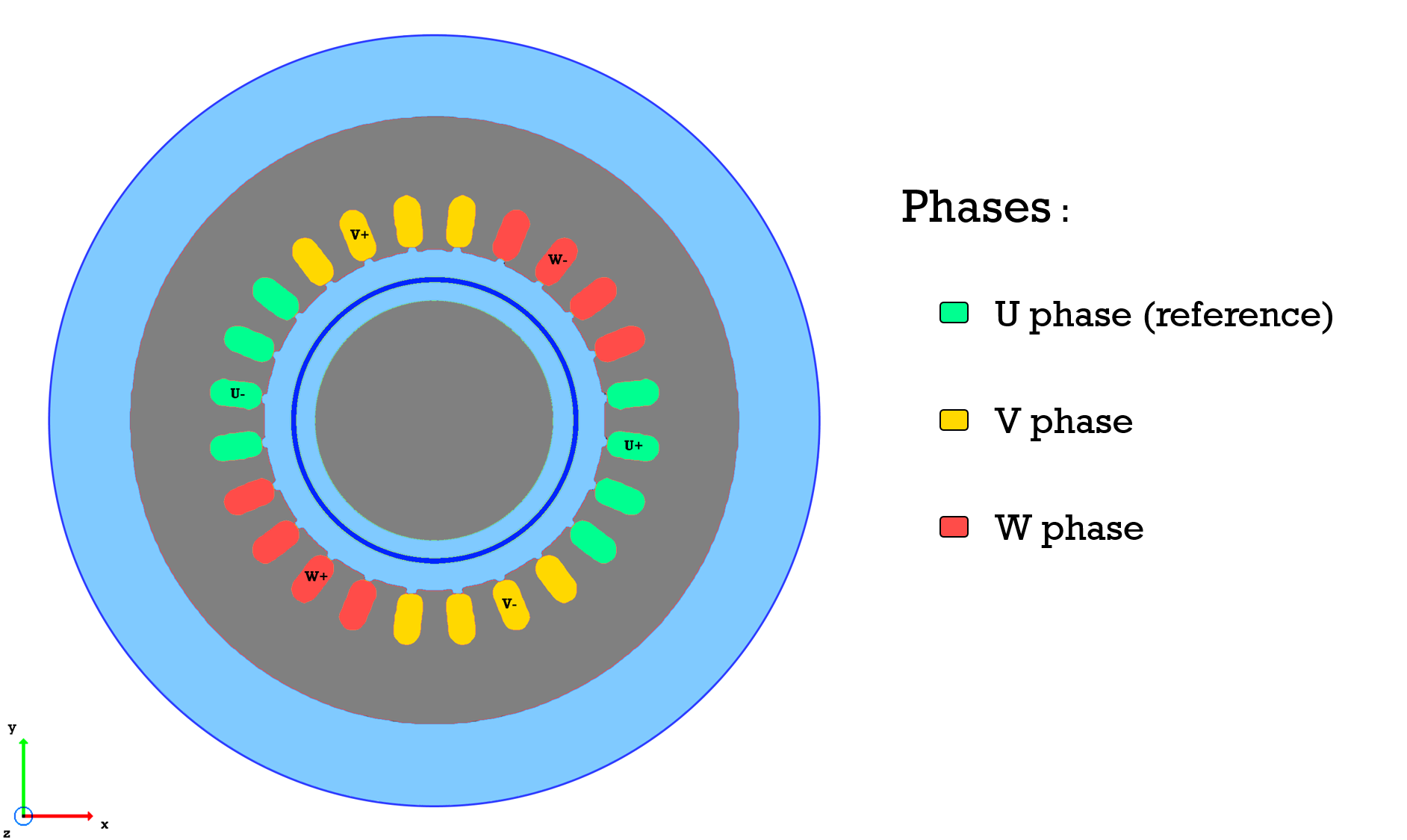}
\caption{Statoric winding distribution and current parameters}
\label{fig:CoilsDistribution}
\end{figure}

The electrical machine geometry and current density distributions are given respectively in Figure \ref{fig:machgeom} and \ref{fig:CoilsDistribution}. The dimensions for the considered machine are given in Table \ref{tab1}.

\begin{table}
\caption{Geometric parameters \label{tab1}}      
\begin{tabular}{ll}
\hline\noalign{\smallskip}
\textbf{Parameter} & \textbf{Value} \\
\noalign{\smallskip}\hline\noalign{\smallskip}
        Slot number         & 24 \\
        Axial length        & 50.0 mm \\
        Outer rotor radius  & 18.5 mm\\
        Inner stator radius & 26.5 mm\\
        Outer stator radius & 47.5 mm\\
        Air gap length      & 8.0  mm \\
\noalign{\smallskip}\hline
\end{tabular}
\end{table}

This machine differs from most conventional SynRM by its large air gap which constrains the statoric winding distribution to only one pair of poles (cf. Figure \ref{fig:CoilsDistribution}).

\begin{table}
\caption{Statoric winding parameters \label{tab2}}
\begin{tabular}{ll}
\hline\noalign{\smallskip}
        \textbf{Parameter} & \textbf{Value} \\
\noalign{\smallskip}\hline\noalign{\smallskip}
        Number of turn $N_s$          & 64 \\
        Winding type                  & Distributed \\
        Connection type               & Star \\
        Resistance ($R_{S,20^{\circ}C}$)     & 7.1 $\Omega$\\
        Voltage $U_{eff}$             & 230 V\\
        Peak intensity $I_{max}$      & 12 A \\
        Number of pole pairs $n_{pp}$ & 1 \\
\noalign{\smallskip}\hline
\end{tabular}
\end{table}

We introduce the relationship between the electrical angle $\theta_{elec}$ and the mechanical angle $\theta$
\begin{equation}
    \theta_{elec} = n_{pp}\theta, \
\end{equation}
with $n_{pp}$ the number of pair of poles, here $n_{pp} = 1$.

We define the three phases as follows: 
\begin{equation}
\label{eq:Coils}
    \begin{cases}
        I_U(\theta) = I_{max} cos\left(n_{pp}\theta + \varphi \right) \\
        I_V(\theta) = I_{max} cos\left(n_{pp}\theta + \varphi - \frac{2\pi}{3} \right) \\
        I_W(\theta) = I_{max} cos\left(n_{pp}\theta + \varphi - \frac{4\pi}{3} \right). 
    \end{cases}
\end{equation}
Here, $\varphi$ is the phase angle.
The computational domain $\Omega$ consists of iron, air, permanent magnet and coils,
\begin{equation}
    \Omega = \Omega_{f} \cup \Omega_{air} \cup \Omega_{mag} \cup \Omega_{c}
\end{equation}
where we further subdivide the ferromagnetic and air subdomains into their rotor and stator parts,
\begin{equation}
    \Omega_f = \Omega_{f, stat} \cup \Omega_{f, rot}, \qquad \Omega_{air} = \Omega_{air, stat} \cup \Omega_{air, rot}.
\end{equation}
Moreover, we subdivide the coil subdomains according to the distribution shown in Figure \ref{fig:CoilsDistribution},
\begin{equation}
    \Omega_c = \Omega_{U^+} \cup \Omega_{U^-} \cup \Omega_{V^+} \cup \Omega_{V^-} \cup \Omega_{W^+} \cup \Omega_{W^-}.
\end{equation}

\subsection{Partial differential equation}
In the two-dimensional magnetostatic setting, the magnetic flux density $\mathbf{ B} = \mbox{curl}((0,0,u)^\top)$ for rotor position $\theta\in [0, 2 \pi]$ can be computed via the solution of the boundary value problem 
\begin{equation}
\begin{aligned}
\label{eq:2dMag}
    \mbox{Find } u \in H^1_0(\Omega): \int_\Omega \nu_\theta(x,|\nabla u|)\nabla u \cdot \nabla v \, dx =\\ \int_{\Omega_{c}} j(\theta) \, v \, dx + \int_{\Omega_{mag}^\theta} R_\theta\begin{bmatrix} -M_y \\ M_x \end{bmatrix} \cdot \nabla v \, dx
\end{aligned}
\end{equation}
for all $v \in H^1_0(\Omega)$, see e.g. \cite{monk_finite_2003}. \\
Here, the magnetic reluctivity is a nonlinear function $\hat \nu$ of the flux density $|\mathbf B| = |\nabla u|$ in the ferromagnetic subdomain and a constant $\nu_0 = 10^7 / (4 \pi)$ elsewhere, i.e.,
\begin{equation} \label{eq_nualpha}
    \nu_{\theta}(x, |\nabla u|) = \begin{cases}
        \hat \nu(|\nabla u|) & x \in \Omega_f^\theta \\
        \nu_0 & x \in \Omega_{air}^\theta \cup \Omega_{c} \cup \Omega_{mag}^\theta
    \end{cases}
\end{equation}
with the rotated domains

\begin{align}
    \Omega_f^\theta =& \Omega_{f,stat} \cup R_\theta \Omega_{f,rot} \\
    \Omega_{air}^\theta =& \Omega_{air,stat} \cup R_\theta \Omega_{air,rot} \\
    \Omega_{mag}^\theta =& R_\theta \Omega_{mag}
\end{align} and $R_\theta$ a rotation matrix around angle $\theta$,
\begin{equation}
    R_\theta = \begin{bmatrix} \mbox{cos}\theta & -\mbox{sin}\theta \\ \mbox{sin}\theta & \mbox{cos}\theta \end{bmatrix} .
\end{equation}
The first term on the right hand side of \eqref{eq:2dMag} represents the impressed current density which is given by
\begin{equation}
\label{eq:current}
    \begin{split}
    j(x,\theta) =& \; \chi_{\Omega_{U^+}}(x) j_U(\theta) +     \chi_{\Omega_{V^+}}(x) j_V(\theta) + \chi_{\Omega_{W^+}}(x) j_W(\theta) \\
    &-\chi_{\Omega_{U^-}}(x) j_U(\theta) - \chi_{\Omega_{V^-}}(x) j_V(\theta) - \chi_{\Omega_{W^-}}(x) j_W(\theta), 
    \end{split}
\end{equation}
where $\chi_A$ denotes the characteristic function of a set $A$,
\begin{equation*}
    \chi_A(x) = \begin{cases}
    1 & x \in A, \\
    0 & \mbox{else}.
    \end{cases}
\end{equation*}
Here the current distribution is defined by :
\begin{equation}
    j_{p}(\theta) = \frac{1}{S_{slot}} N_s I_p(\theta), \quad p \in \{U,V,W\}
\end{equation}
with $S_{slot}$ the cross-sectional area of one coil, $N_s$ the number of turns per coil and $I_U$, $I_V$, $I_W$ as defined in \eqref{eq:Coils}.
The second term on the right hand side of \eqref{eq:2dMag} represents the magnetization $\mathbf M = (M_x, M_y)^\top$ coming from permanent magnets which will be added in the course of the multi-material optimization procedure.

In the following, we will denote by $u_{\theta}$ the solution to the state equation \eqref{eq:2dMag}.\\
We present here after the properties of interest of the materials (air,ferromagnetic,magnet) used in the machine. 

\begin{table}[H]
\caption{Material properties \label{tab3}}
\begin{tabular}{lll}
\hline\noalign{\smallskip}
        \textbf{Material} & \textbf{Reluctivity [$m.H^{-1}$]} & \textbf{Magnetization [$A.m^{-1}$]}  \\
\noalign{\smallskip}\hline\noalign{\smallskip}
        Air & $\nu_0$ & 0 \\
        Copper & $\nu_0$ & 0 \\
        Ferromagnetic & $\hat \nu(|\vec{B}|)$ & 0 \\
        Magnet & $\nu_0$ & $M_{max}$ \\
\noalign{\smallskip}\hline
\end{tabular}
\end{table}

\begin{figure}
     \begin{subfigure}[b]{0.4\textwidth}
         \includegraphics[width=6 cm]{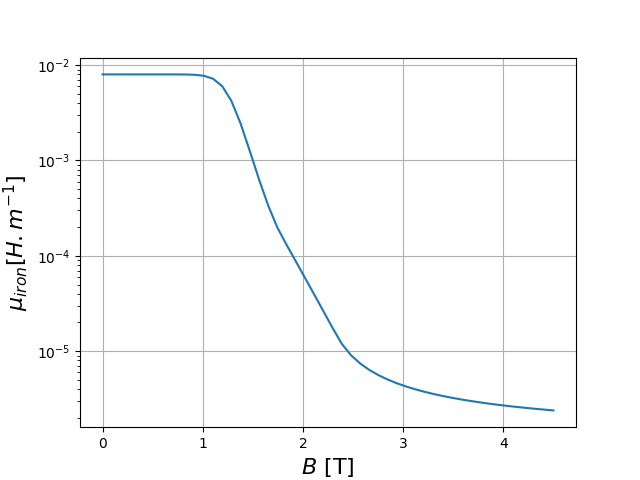}
         \caption{ Magnetic permeability}     
     \end{subfigure}
     \hfill
     \begin{subfigure}[b]{0.4\textwidth}
        \includegraphics[width=6cm]{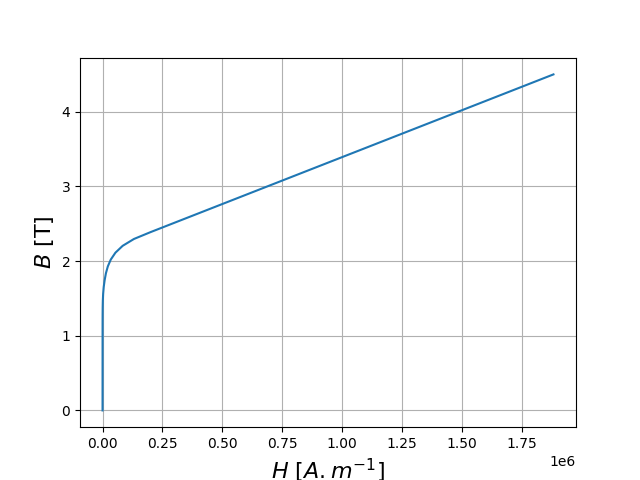}
         \caption{Marrocco BH-curve}
         \label{fig:five over x}
     \end{subfigure}
\caption{Non-linear characteristic for the ferromagnetic material.}     
        \label{fig:relucBH}
\end{figure}

The maximum norm of the magnetization vector was chosen as $M_{max} = 2.33 \cdot 10^5 A.m^{-2}$ to fit data from \cite{nunes_modeling_2020} on ferrite magnets.
The reluctivity of the magnets and of the copper coils is assimilated to the one of air to simplify further material interpolation and avoid complex schemes like the ones found in \cite{zuo_multi-material_2017}. The non-linear behaviour of the ferromagnetic material is modelled with a Marrocco's BH curve approximation \cite{marrocco_analyse_1977}.
\begin{equation}
    \hat  \nu (|\vec{B}|) =
    \begin{cases}
    \nu_0 (\varepsilon + \frac{(c-\varepsilon)|\vec{B}|^{2\alpha}}{\tau + |\vec{B}|^{2\alpha}} `& \mbox{if  }|\vec{B}|\leq B_{max}, \\
    \nu_0 \left( 1 - \frac{M_s}{|\vec{B}|} \right) & \mbox{else if  } |\vec{B}|> B_s, \\
    exp\left(\frac{\gamma \left(|\vec{B}| - \beta \right)}{|\vec{B}|}\right) & \mbox{otherwise},
    \end{cases}
\end{equation}
where $B_s = \beta + \frac{log \left( \frac{\nu_0}{\gamma} \right)}{\gamma}$ and $M_s = B_s + \frac{1}{\gamma}$ and the coefficient of the Marrocco curve in Figure\ref{fig:relucBH} are defined in the Table \ref{tab4}.

\begin{table}[H]
\caption{Marrocco curve coefficient for the ferromagnetic material\label{tab4}}
\begin{tabular}{lll}
\hline\noalign{\smallskip}
        \textbf{Parameter} & &\textbf{Value}  \\
\noalign{\smallskip}\hline\noalign{\smallskip}
        \textbf{$\alpha$} & & 6.84 \\
         \textbf{$\beta$} & & -1.30$\cdot10^{-1}$ \\
         \textbf{$\gamma$} & & 4.86 \\
         \textbf{$\varepsilon$} & & 1.57 $\cdot10^{-4}$ \\
          \textbf{$\tau$} & & 4.14 $\cdot10^{3}$ \\
         \textbf{c}  & & 1.90$\cdot10^{-2}$ \\
         \textbf{$B_{max}$} & & 1.80 (T)\\
\noalign{\smallskip}\hline
\end{tabular}
\end{table}

\subsection{Torque computation method}
For computing the torque, we chose a method based on Maxwell's stress tensor, Arkkio's method \cite{sadowski_finite_1992}. 
While co-energy torque computation methods were proven to be more precise and less prone to error, they are more costly in terms of computation time and not fit for optimization.
Using Arkkio's method, the torque can be computed as

\begin{equation}
    T =  \frac{L_z \nu_0}{r_s-r_r} \int_{S} \sqrt{x^2+y^2} B_r B_\phi dS
\end{equation}
where $B_r$ and $B_\phi$ denote the radial and tangential magnetic induction, respectively, $L_z$ denotes the length of the machine in $z$-direction and $S$ denotes the surface between radii $r_s$ and $r_r$ in the air gap (with $r_s>r_r$).

In the setting of two-dimensional magnetostatics, the torque for the rotor position given by angle $\theta$ thus amounts to 
\begin{equation}
    T_{\theta} = T(u_{\theta}) =  \frac{L_z \nu_0}{r_s-r_r} \int_{S} Q\nabla u_{\theta} \cdot \nabla u_{\theta} \, dS
\end{equation}
\begin{equation}
    \mbox{with } Q(x,y) = 
    \begin{bmatrix} \frac{xy}{\sqrt{x^2+y^2}} && \frac{y^2-x^2}{2\sqrt{x^2+y^2}} \\   \frac{y^2-x^2}{2\sqrt{x^2+y^2}} && \frac{-xy}{\sqrt{x^2+y^2}}
    \end{bmatrix} \in \mathbb R^{2 \times 2}.
\end{equation}

Determining the average torque by means of its instantaneous values can be very expensive. It is shown in \cite{bianchi_mmf_2007} that a good approximation to the average torque can be obtained when evaluating the torque for only suitably chosen rotor positions,
\begin{equation} \label{eq:torque4point}
    \bar{T} = \frac{1}{4} \left(T_{0} + T_{\frac{\pi}{12}} + T_{\frac{\pi}{6}} + T_{\frac{\pi}{4}} \right).
\end{equation}
We compared the average torque obtained by evaluation at 500 equally distributed rotor positions between $0$ and $2\pi$ with the value obtained by the four-point formula \eqref{eq:torque4point}.

When the torque value is not equal to zero, the error found was to be lower than 0.1\% as expected and described in literature \cite{Akiki8507245}. This is solved beyond the first iteration.
\begin{table}
\caption{Four static positions method error \label{tab5}}
\begin{tabular}{lllll}
\hline\noalign{\smallskip}
        Design & \textbf{$\bar{T} [N.m]$} & \textbf{$\bar{T} [N.m]$} & \textbf{Error}  \\
        & \textbf{(500 points)}  & \textbf{(4points)}&  \textbf{[\%]}\\
 \noalign{\smallskip}\hline\noalign{\smallskip}
        Unbiased starting point* &  1.5790 $\cdot 10^{-6}$ & 5.7232  $\cdot 10^{-6}$ & 262.4 \\
        Final design Table \ref{IronTable} & 1.1123 &  1.1129 & 0.048 \\
        Final design Table  \ref{IronMagnetNeutralTable} & 1.4513 & 1.4516 & 0.027 \\
\noalign{\smallskip}\hline
\end{tabular}
{\raggedright * in this design $\rho_\nu =0.5,\rho_{M_x} =0.5,\rho_{M_y} =0.5$ everywhere in the rotor \par}
\end{table}

\section{Optimization problem}
In this section, we define our optimization problem and reformulate the forward problem to fit the density-based topology optimization approach. Our goal is to maximize the average torque computed via \eqref{eq:torque4point},
\begin{equation}
    (P_1):
    \begin{cases}
    \text{maximize } \bar{T} =\frac{1}{4} \left(T(u_{0}) + T(u_{\frac{\pi}{12})} + T(u_{\frac{\pi}{6}}) + T(u_{\frac{\pi}{4}}) \right) \\
    \mbox{s.t. $u_{\theta}$ is a solution of \eqref{eq:2dMag} for }\theta \in \{0,\frac{\pi}{12},\frac{\pi}{6},\frac{\pi}{4}\}
    \end{cases}
\end{equation}
This is achieved by finding the optimal material distribution consisting of ferromagnetic material, air and permanent magnets on the one hand, and the optimal magnetization direction of the permanent magnets on the other hand. Moreover, we will incorporate a bound on the maximum allowed permanent magnet volume.
\subsection{Density based topology optimization}
Let us reformulate the forward problem \eqref{eq:2dMag}, introducing the three density variables respectively for the ferromagnetic material and the two components of the permanent magnets magnetization, $\rho_\nu,\rho_{M_x},\rho_{M_y}$ defined in $\Omega^{\theta}_{rot} = R_\theta ( \Omega_{f,rot} \cup \Omega_{air,rot} \cup \Omega_{mag}$). Moreover we introduce the rotated design variables
\begin{align*}
    \rho_\nu^\theta(x,y) =& \rho_\nu( R_\theta((x,y)^\top)) \\
    \rho_{M_x}^\theta(x,y) =& \rho_{M_x}( R_\theta((x,y)^\top)) \\
    \rho_{M_y}^\theta(x,y) =& \rho_{M_y}( R_\theta((x,y)^\top)) 
\end{align*}
which represent the design given by $\rho_\nu$, $\rho_{M_x}$, $\rho_{M_y}$ after rotation, and the vector of design variables $\textbf{X} := \begin{bmatrix}\rho_\nu, & \rho_{M_x}, & \rho_{M_y} \end{bmatrix}^\top $.\\

Given two interpolation functions
\begin{equation}
f_\nu: [0,1]  \rightarrow [0,1], \qquad f_M: [0,1]  \rightarrow [0,1],
\label{eq:Frho}
\end{equation}
we define the operator 
\begin{equation}
\begin{aligned}
\label{eq:OpPDE_density}
    K_{\theta}: (\mathbf{X},u,v) \mapsto &
    \int_{\Omega}\nu(\rho_{\nu}^\theta,|\nabla u|) \nabla u \cdot \nabla v \\
    & - \int_{\Omega_{rot}} f_{\nu}(1-\rho_{\nu}^\theta) \frac{M_{max}f_{M}(|\vec{M}^\theta|)}{|\vec{M}^\theta|}
    R_{\theta} \begin{bmatrix} -M_y^\theta \\ M_x^\theta \end{bmatrix} \cdot \nabla v,&
\end{aligned}
\end{equation}
with the reluctivity function 
\begin{equation}
\nu(\rho_{\nu}^\theta,|\nabla u|) = 
\begin{cases}
\hat \nu(|\nabla u|) & \text{ in }\Omega_{f,stat} \\
\nu_0                & \text{ in }\Omega_{c} \cup \Omega_{air,stat}  \\
\nu_0 + f_{\nu}(\rho_{\nu}^\theta)(\hat \nu(|\nabla u|) - \nu_0) &
\text{ in } \Omega_{rot}
\end{cases}
\end{equation}
and with the components of the magnetization vector $\vec{M}^\theta = (M_x^\theta, M_y^\theta)$ given in dependence of the two rotated density variables $\rho_{M_x}^\theta$, $\rho_{M_y}^\theta$,
\begin{equation} \label{eq:def_fsd_tilde}
 (M_x^\theta,M_y^\theta) = \tilde f_{sd}(\rho_{M_x}^\theta, \rho_{M_y}^\theta) 
\end{equation}
for a mapping $\tilde f_{sd}$ which will be discussed later on.
Hence, the state equation \eqref{eq:2dMag} can be reformulated into 
\begin{equation}
\begin{aligned}
\label{eq:PDE_density}
&\text{Find } u \in H^1_0(\Omega): \, \\
&K_{\theta}(u,v,\textbf{X}) = 
\int_{\Omega_{c}}j(\theta) \, v \, dx, \text{ for all } v \in H^1_0(\Omega).
\end{aligned}
\end{equation}

The optimization problem ($P_1$) can then be reformulated into
\begin{equation}
    (P_2):
    \begin{cases}
    \underset{\textbf{X}}{\text{maximize }}  \bar{T} = \frac{1}{4} \sum_{\theta \in \{ 0,\frac{\pi}{12},\frac{\pi}{6},\frac{\pi}{4} \}}T(u_{\theta}) \\
    \text{s.t. $u_{\theta}$ is a solution of the \eqref{eq:PDE_density}} \quad \theta \in \{ 0,\frac{\pi}{12},\frac{\pi}{6},\frac{\pi}{4} \}
    \end{cases}
\end{equation}


\subsubsection{Material interpolation schemes}
In density based topology optimization, the quality of the final solution is dependant on the choice of interpolation functions (equation \eqref{eq:Frho}). We present here two existing schemes and a novel one based on properties of the topological derivative. \\
The polynomial interpolation scheme 
\begin{equation}
\label{eq:poly}
    f_n(\rho) = \rho^n \quad n>0,
\end{equation}
also referred to as SIMP (Solid Isotropic Material with Penalization), is the most used material interpolation scheme for topology optimization and allows for easy penalization of intermediate materials.
However, it presents some symmetry issues and favors low $\rho$ associated material in the final design. In \cite{sanogo_topology_2018}, the authors compared this scheme to other schemes and concluded that the final design was not as good as many other proposed ones.
\begin{figure}[H]
\includegraphics[width=8 cm]{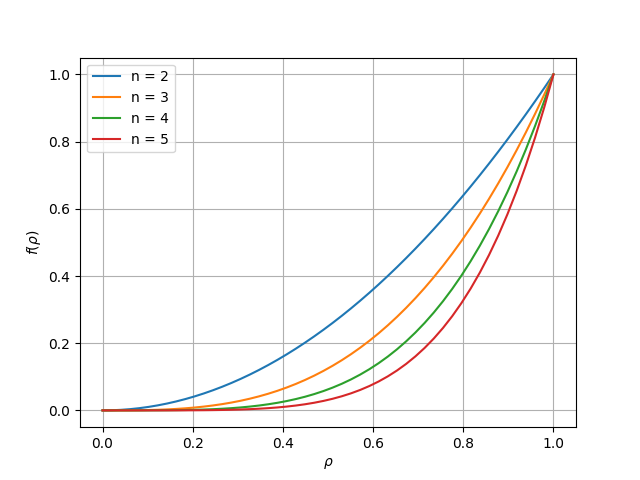}
\caption{SIMP Polynomial interpolation scheme}
\label{fig:PolyInterpol}
\end{figure} 

\begin{figure}
\centering
\includegraphics[width=8 cm]{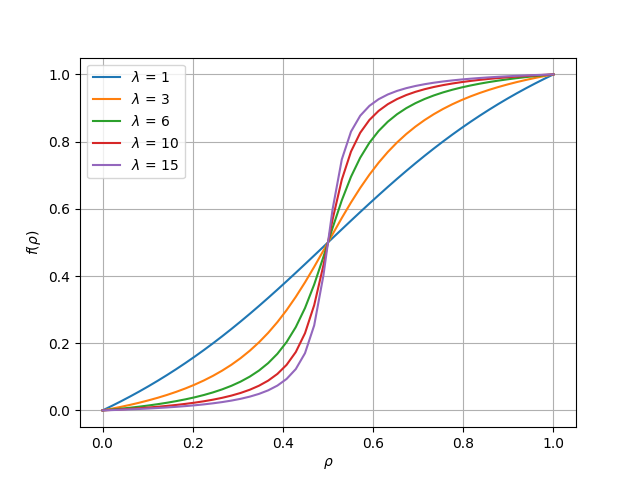}
\caption{D. Luk\`{a}\v{s}'s interpolation scheme }
\label{fig:LukasInterpol}
\end{figure}

\begin{figure}
\includegraphics[width=8 cm]{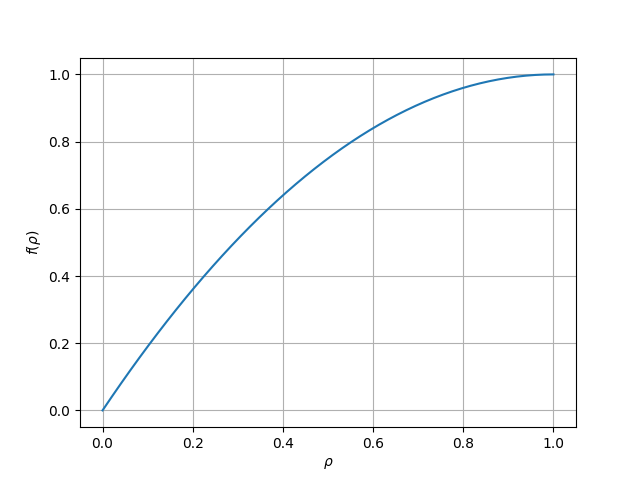}
\caption{Topological derivative inspired interpolation scheme }
\label{fig:GanglInterpol}
\end{figure}

To solve symmetry issues introduced by the classical polynomial interpolation,
D. Luk\`{a}\v{s} introduced a new scheme in \cite{lukas_integration_2006}: 
\begin{equation}
    f_\lambda(\rho) = \frac{1}{2}\left(1 + \frac{1}{\arctan(\lambda)} \arctan(\lambda(2\rho-1))\right), \quad \lambda>0.
\end{equation}
In this equation the particular invariant point $\rho=0.5$ does not promote intermediate materials, grey material depends on $\lambda$ values (cf. Figure \ref{fig:LukasInterpol}). \\
High $\lambda$ values permit to penalize intermediate materials but can lead to a poor convergence of the algorithm. A parameter study for $\lambda$ led us to choose $\lambda=5$. This interpolation method is chosen for the norm of the magnetization vector ($f_M$ in \eqref{eq:PDE_density}).
Finally, we propose a new interpolation scheme as given in Figure \ref{fig:GanglInterpol},which is inspired by the topological derivative as done in \cite{AmstutzDapognyFerrer2018}, see also the the SIMP-All method for linear elasticity \cite{FerrerSIMPAll}. Here, we seek to design a material interpolation function whose derivative with respect to the density variable $\rho$ coincides with the topological derivative of the problem at $\rho =0$ and $\rho = 1$. When interpolating between two linear materials with reluctivity values $\nu_0$ and $\nu_1$, the conditions for the material interpolation function $f$ according to \cite{AmstutzDapognyFerrer2018} would read
\begin{equation} \label{eq_conditions_f}
    \begin{cases}
    f(0) = 0, \\
    f(1) = 1, \\
    f'(0) = \frac{2\nu_0}{\nu_0 + \nu_{1}}, \\
    f'(1) = \frac{2\nu_{1}}{\nu_0 + \nu_{1}}. 
    \end{cases}
\end{equation}
Due to the involved formula of the topological derivative for nonlinear magnetostatics \cite{AmstutzGangl2019}, a mathematically rigorous extension of this method to the nonlinear setting is not straightforward. However, inspired by the particular behaviour of the Marrocco BH-curve where the magnetic reluctivity is almost constant for low flux density values, see Figure \ref{fig:relucBH}, we simply use this idea for that constant reluctivity value $\nu_1:= \nu_0 \varepsilon \approx 124.94$. Using cubic Hermite interpolation for the conditions \eqref{eq_conditions_f}, we obtain the polynomial
\begin{equation} \label{eq_mif_td}
    f(\rho) = \frac{2\nu_0}{\nu_0 + \nu_{1}} \rho 
    - \frac{\nu_0-\nu_{1}}{\nu_0 + \nu_{1}} \rho^2.
\end{equation}
Note that the term of order 3 happens to vanish.

 \begin{figure}
     \begin{subfigure}[b]{0.4\textwidth}
         \includegraphics[width=7cm]{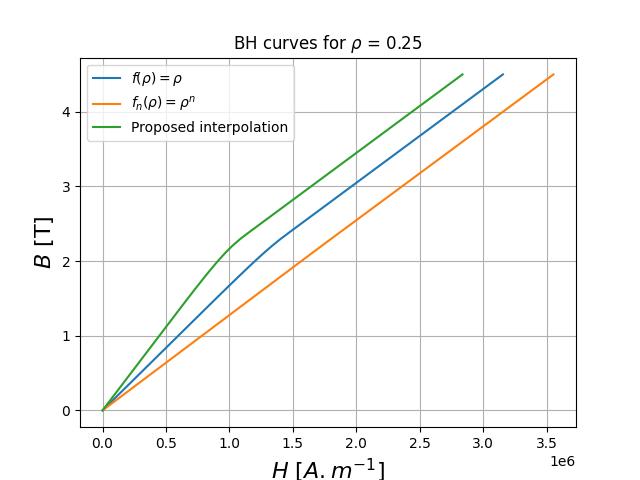}
         \caption{Interpolated BH curve for $\rho=0.25$}
                  \label{fig:BH curve025}
     \end{subfigure}
     
     \begin{subfigure}[b]{0.4\textwidth}
        \includegraphics[width=7cm]{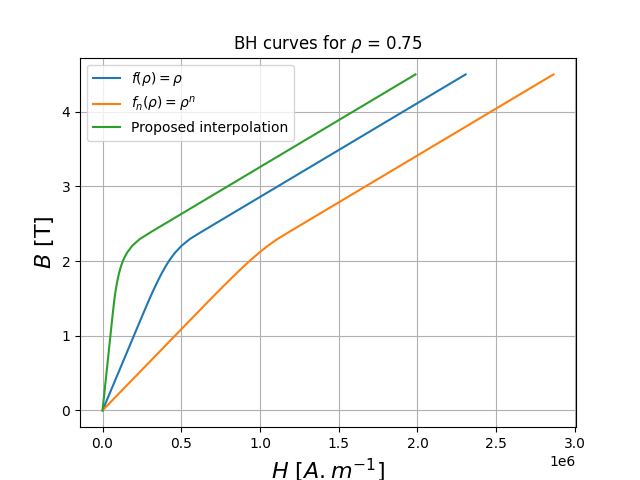}
         \caption{Interpolated BH curve for $\rho=0.75$}
         \label{fig:BH curve075}
     \end{subfigure}
     
     \begin{subfigure}[b]{0.4\textwidth}
        \includegraphics[width=7cm]{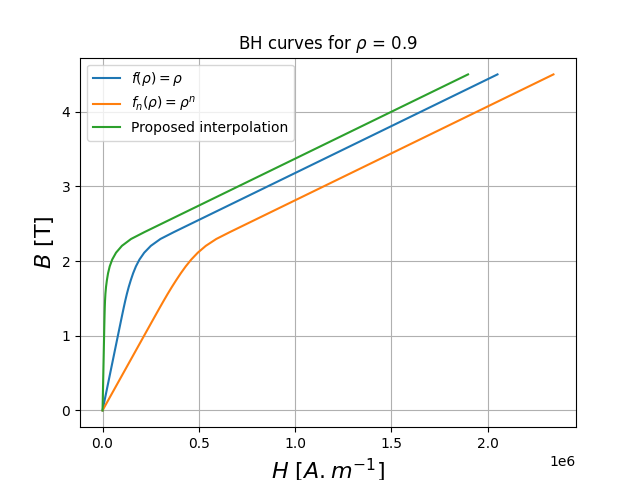}
         \caption{Interpolated BH curve for $\rho=0.9$}
         \label{fig:BH curve09}
     \end{subfigure}
     
\caption{Interpolated BH curves using material interpolation scheme \eqref{eq_mif_td}.}   

        \label{fig:BHinterpol}
        
\end{figure}

\subsubsection{Magnetization vector transform}
We deal with two magnetization density variables $\rho_{M_x}$ and $\rho_{M_y}$ in order to represent the magnetization direction $(M_x,M_y)$. One way of relating these quantities to each other would be to have $\rho_{M_x}$ represent the first and $\rho_{M_y}$ the second coordinate, resulting in a representation in Cartesian coordinates, which was also considered in \cite{wang_topology_2005}. In this case, however, some magnetization directions exhibit higher maximum magnetization than others, e.g. $\rho_{M_x} = \rho_{M_y}=1$ would correspond to $| (M_x, M_y)^\top | = \sqrt{2}$ whereas for the magnetization direction pointing to the right $\rho_{M_x} =1$, $ \rho_{M_y}=0.5$ would yield a maximum magnetization of $| (M_x, M_y)^\top | = 1$, thus making the maximum magnetization angle dependent.
As an alternative, one could use polar coordinates and represent the magnetization direction by just one periodic density function. In this case, however, the ambiguity of angular values causes problems in the gradient computation.
We define a means of solving these issues without resorting to polar coordinates. 
We decide on two density variables $\rho_{M_x}, \rho_{M_y}$ with values in $[0,1]$, but map them onto a disk, thereby avoiding an angle dependent maximum magnetization value.

\begin{figure}
\centering
\includegraphics[width=8 cm]{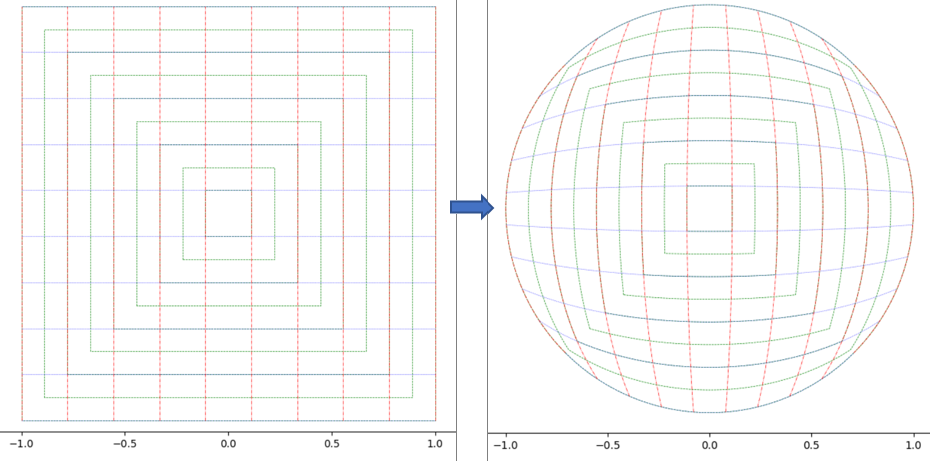}
\caption{Square to disk transform for the Magnetization vector coordinates}
     \label{fig:/Mapping}
\end{figure}
There are several mappings that approximately realize such a square-to-disk transformation. In order to preserve the angle,as much as possible, without being too computationally heavy, the elliptic grid mapping.
\begin{equation}
\label{eq:ellipticSC}
f_{sd}(x,y) = 
\begin{cases}
x\sqrt{1-\frac{y^2}{2}} \\
y\sqrt{1-\frac{x^2}{2}}
\end{cases} \text{with } (x,y) \in [-1,1]^2. 
\end{equation}

was chosen as a good compromise out of the methods described in \cite{fong_analytical_2019}.
While the associated inverse transformation $f_{ds} = f_{sd}^{-1}$ given by

\begin{equation}
\label{eq:ellipticCS}
f_{ds}(u,v) = 
\begin{cases}
\frac{1}{2}\left(\sqrt{2+u^2-v^2+2\sqrt{2}u} - \sqrt{2+u^2-v^2-2\sqrt{2}u}\right) \\
\frac{1}{2}\left(\sqrt{2-u^2+v^2+2\sqrt{2}v} - \sqrt{2-u^2+v^2-2\sqrt{2}v}\right)
\end{cases}
\end{equation}

is computationally more costly, it is only used once per iteration and in post-processing.
The mapping between the magnetization density variables $\rho_{M_x}, \rho_{M_y}$ and the magnetization vector $\vec{M} = (M_x, M_y)$ \eqref{eq:def_fsd_tilde} is then given by
\begin{align} \label{eq:def_fsd_tilde2}
    \tilde f_{sd}(\rho_{M_x}^\theta, \rho_{M_y}^\theta) = f_{sd}(2(\rho_{M_x}^\theta-0.5),2(\rho_{M_y}^\theta-0.5))
\end{align}

\subsection{Incorporation of volume constraints}
Constraints on iron and magnets volume are added to avoid having structures with disproportionate volumes of material. Hence we add new constraints with the operator
\begin{equation}
    I_{vol}: \rho \mapsto \frac{1}{V_{\Omega_{rot}}}\int_{\Omega_{rot}} \rho(x) \, dx
\end{equation}
representing the volume fraction inside the rotor domain $\Omega_{rot}$ of a material given by a density function $\rho$. Here, $V_{rot} = \int_{\Omega_{rot}} 1 \; dx$ denotes the total area of $\Omega_{rot}$.
Based on ($P_2$) we define a new constrained optimization problem
\begin{equation}
    (P_3): 
    \begin{cases}
    \text{maximize }  \bar{T} = \frac{1}{4} \sum_{\theta \in \{ 0,\frac{\pi}{12},\frac{\pi}{6},\frac{\pi}{4} \}}T(u_{\theta}) \\
    \text{s.t. $u_{\theta}$ is a solution of $\eqref{eq:PDE_density}, \; \theta \in \{ 0,\frac{\pi}{12},\frac{\pi}{6},\frac{\pi}{4} \} $} 
    \\
    \quad \; \text{ $I_{vol}(\rho_{\nu})\leq f_{v,f}$} 
    \\
    \quad \; \text{ $I_{vol}((1-\rho_{\nu}) \, |\vec{M}|)\leq f_{v,mag}$}
    \end{cases}
\end{equation}
with given upper bounds on the allowed ferromagnetic and permanent magnet material $f_{v,f}$, $f_{v,mag} \in [0,1]$, respectively.
We reformulate the inequality constraints of ($P_3$) using the augmented Lagrangian framework as described in \cite{nocedal_penalty_1999},
\begin{equation}
\label{eq:AugLag}
    (P_4): 
   \begin{cases}
    \text{minimize }
    L(\mathbf{X}, \underline u) = 
    -\bar{T}(\underline u)
    + \psi(h_{v,f}(\mathbf{X}),\gamma_f,\mu) \\
          \quad \quad \quad \quad \quad \quad \quad \quad
    + \psi(h_{v,mag}(\mathbf{X}),\gamma_{mag},\mu)&\\
    \text{s.t. $u_{\theta}$ is a solution of \eqref{eq:PDE_density}},  
    \theta \in \{ 0,\frac{\pi}{12},\frac{\pi}{6},\frac{\pi}{4} \}  \\
    \end{cases}
\end{equation}
With the state vector
\begin{align} \label{eq_state_vector}
    \underline u := (u_0, u_{\frac\pi{12}}, u_{\frac\pi6}, u_{\frac\pi4}).
\end{align}
Here, 
\begin{equation}
    \begin{cases}
    h_{v,f}(\mathbf{X}) = f_{v,f}-I_{vol}(\rho_\nu), \\
    h_{v,mag}(\mathbf{X}) = f_{v,mag}-I_{vol}((1-\rho_{\nu}) \,|\vec{M}|), \\
    \end{cases}
\end{equation}
with the scalar function
\begin{equation}
    \psi(t,\sigma,\mu) = 
    \begin{cases}
    -\sigma t + \frac{1}{2\mu}t^2 & \text{if } t-\mu \sigma \leq 0, \\
    -\frac{\mu}{2}\sigma^2 & \mbox{otherwise}.
    \end{cases}
\end{equation}
The positive scalar multipliers $\gamma_{f},\gamma_{mag},\mu$ are updated according to the LANCELOT-Method of Multipliers presented in \cite{nocedal_penalty_1999}.

\subsection{Adjoint Method}

To solve the optimization problem $(P_4)$ as formulated in \eqref{eq:AugLag} by a gradient descent algorithm, we introduce the Lagrangian for the PDE-constrained problem \eqref{eq:AugLag}
\begin{equation}
\label{eq:unconstrLag}
    \mathcal L(\mathbf X, \underline u, \underline w) = L (\mathbf X, \underline u) + \sum_{\theta \in \{ 0,\frac{\pi}{12},\frac{\pi}{6},\frac{\pi}{4} \}} K_\theta(\mathbf X,  u_\theta, w_\theta) - \int_{\Omega_{stat}} j(\theta) w_\theta,
\end{equation}
where $\underline w = (w_0, w_{\frac\pi{12}}, w_{\frac\pi6}, w_{\frac\pi4})$ is a vector of Lagrange multipliers. The adjoint states $\lambda_\theta$ corresponding to problem \eqref{eq:AugLag} for different rotor positions $\theta$ are defined by
$\frac{\partial \mathcal L}{\partial u_\theta}(\mathbf X, u_\theta, \lambda_\theta) = 0$, i.e.
$\lambda_\theta$ is the solution to

\begin{equation}
\label{eq:adjoint}
\mbox{Find } \lambda_\theta \in H_0^1(\Omega): \frac{\partial K_{\theta}}{\partial u}(\mathbf X, u_\theta, \lambda_\theta)(v) = - \frac14 \frac{\partial T}{\partial u} (u_\theta)(v) \quad \\ \mbox{for all } v \in H_0^1(\Omega).
\end{equation}

Here, the left and right hand side, respectively, involve the operators
\begin{equation}
\begin{aligned}
    \frac{d}{du} K_\theta(\mathbf X, u, \lambda)(v) =&  \frac{d}{du} \left( \int_{\Omega}\nu(\rho_{\nu}^\theta,|\nabla u|) \nabla u \cdot \nabla \lambda \right) (v) \\
    =& \int_{\Omega}\nu(\rho_{\nu}^\theta,|\nabla u|) \nabla v \cdot \nabla \lambda   +   \int_{\Omega} \frac{d}{du} \nu(\rho_{\nu}^\theta,|\nabla u|)(v) \nabla u \cdot \nabla \lambda \\
    =& \int_{\Omega}\nu(\rho_{\nu}^\theta,|\nabla u|) \nabla v \cdot \nabla \lambda   +   \int_{\Omega} f_\nu(\rho_\nu^\theta) \frac{\hat \nu'(|\nabla u|)}{|\nabla u|} (\nabla u \cdot \nabla v) (\nabla u \cdot \nabla \lambda) \\
    \frac{d}{du}T(u)(v) =& \; 2 \frac{L_z \nu_0}{r_s - r_r} \int_S Q \nabla u \cdot \nabla v \; dS.
\end{aligned}   
\end{equation}
Similarly to \eqref{eq_state_vector}, we introduce the adjoint vector
\begin{equation}
    \underline \lambda := (\lambda_0, \lambda_{\frac\pi{12}}, \lambda_{\frac\pi6}, \lambda_{\frac\pi4}).
\end{equation}
Hence, the sensitivity associated with \eqref{eq:unconstrLag} amounts to
\begin{equation}
\label{eq:diffLag}
    \begin{split}
     \frac{\partial \mathcal{L}}{\partial \mathbf{X}} (\mathbf X, \underline u, \underline \lambda) = & 
     \sum_{\theta \in \{0,\frac{\pi}{12},\frac{\pi}{6},\frac{\pi}{4}\}} \frac{\partial K_\theta}{\partial \mathbf{X}}(\mathbf X, u_\theta,\lambda_\theta)  \\
     & + 
     \frac{\partial \psi(h_{v,f}(\mathbf{X}),\gamma_f,\mu)}{\partial \mathbf{X}} +
     \frac{\partial \psi(h_{v,mag}(\mathbf{X}),\gamma_{mag},\mu)}{\partial \mathbf{X}}
    \end{split}
\end{equation}
where we used that $\frac{\partial T}{\partial \mathbf X} = 0$ and $ \frac{\partial j}{\partial \mathbf X} = 0$.


\subsection{Update method}
Several methods exist to consider the movement in electrical motors, such as the Moving Band (MB) technique described in \cite{Davat1985}. Even when using high order elements in the MB, this method remains less accurate than the mortar element method \cite{Antunes2005}. A variant of the mortar method, the Nitsche method, is chosen to take into account the rotation \cite{hollaus_nitsche-type_2010}. The operator $K_\theta$ will be replaced by $K_\theta^{NM}$ and $\mathcal{L}$ by $\tilde{\mathcal{L}}$ to fit the new formulation (detailed in Appendix.\ref{appendix}).

We introduce a triangular mesh with a total of $N$ elements, with $N_{rot}$ many elements inside the rotor domain $\Omega_{rot}$. We use piecewise linear and globally continuous finite element basis functions to solve the state and adjoint equations and we represent the density variables $\rho_{\nu}, \rho_{M_x}, \rho_{M_y}$ as piecewise constant functions on the mesh corresponding to $\Omega_{rot}$. Thus, these density variables can be represented by a vector of dimension $N_{rot}$ consisting of the values of the discrete functions in each element. We will use the same notation $\rho_{\nu}, \rho_{M_x}, \rho_{M_y}$ for the vectors representing the discretized density variables.

To comply with the bounds of the density variables, we introduce the projected gradient as defined in \cite[p. 520]{nocedal_penalty_1999} 
\begin{equation}
\label{eq:projectedGrad}
\mathcal{G} := 
\begin{bmatrix} 
\mathbf{P}_{\rho_\nu,[0,1]} \left(\frac{\partial \tilde{\mathcal{L}}}{\partial \rho_\nu} \right) & 
\mathbf{P}_{\rho_{M_x},[0,1]} \left(\frac{\partial \tilde{\mathcal{L}}}{\partial \rho_{M_x}} \right) & 
\mathbf{P}_{\rho_{M_y},[0,1]} \left(\frac{\partial \tilde{\mathcal{L}}}{\partial \rho_{M_y}} \right) 
\end{bmatrix}^\top,
\end{equation}
with the projection operator $\mathbf{P}_{\rho,[a,b]} : \mathbb R^{N_{rot}} \rightarrow \mathbb R^{N_{rot}}$ defined by
\begin{equation}
    \left(\mathbf{P}_{\rho,[a,b]}(X) \right)_i =
    \begin{cases}
    X_i  & \mbox{if } \rho_i \in ]a,b[, \\
    min(0,X_i)  & \mbox{if } \rho_i = a ,\\
    max(0,X_i)  & \mbox{if } \rho_i = b, \\
    \end{cases}
\end{equation}
for $i = 1,...,N_{rot}$.
This in itself allows for finer geometry to emerge by amplifying the relative importance of the gradient where change in the density function can be made.\\
From the sensitivity we derive the update equation at iteration n
\begin{equation}
    \mathbf{X}_i^{n+1} = Q_{[0,1]} \left( \mathbf{X}_i^{n} - s\frac{\mathcal{G}_i}{|\mathcal{G}_i|} \right) , \quad \mbox{for }i = 1,...,N_{rot}
\end{equation}
with the projection operator $Q_{[a,b]}: \mathbb R^{3 N_{rot}} \rightarrow \mathbb R^{3 N_{rot}}$ defined by
\begin{equation}
    \left(Q_{[a,b]}(\mathbf v) \right)_i = max(a,min(b,v_i))
\end{equation}
for a vector $\mathbf v = (v_1, \dots, v_{N_{rot}})^\top \in \mathbb R^{3 N_{rot}}$ and $i=1, \dots,N_{rot}$, to enforce the bounds of the density variables.
Here, $s$ denotes the step size which is chosen in such a way that a descent of the augmented Lagrangian is obtained,
\begin{equation}
    L(\mathbf X^{n+1}, \underline u^{n+1}) < L(\mathbf X^{n}, \underline u^n)
\end{equation}
with $\underline u^n$ is the vector of states for the design represented by $\mathbf X^n$.

\subsection{Filtering and projection method}

In density-based topology optimization, checkerboard patterns and small isolated elements of one material are avoided using filtering methods at each step of gradient descent. While this filtering has a regularizing effect on the density variables, it may introduce more gray areas. Therefore, in the next step, so-called projection methods are applied in order to get a more defined shape. This combination allows achieving smoother and more defined boundaries between the materials in the final design.

In our approach, in the first step we perform density filtering by solving the PDE
\begin{equation} \label{eq_filterHelm}
    -r^2\nabla^2 \rho + \rho = \rho_{ref} \quad r>0.
\end{equation}
with $\rho_{ref}$ the given density function, which is commonly referred to as Helmholtz filtering \cite{lazarov_filters_2011}.
This mesh-independency filter modifies the sensitivity by averaging on the neighbor cells. The parameter $r$ is the radius of influence of the filtering and is defined in our case as a factor $\delta$ of the minimum representative mesh-element length $h$, i.e. $r = \delta \; h$. \\

For the projection step, we choose the function proposed in \cite{wang_projection_2011},
\begin{equation} \label{eq_projection}
\begin{aligned}
    &f_{\rho_{cut},\beta}(\rho) = \frac{\tanh(\beta(\rho-\rho_{cut})) + \tanh(\beta \rho_{cut})}
    {\tanh(\beta(1-\rho_{cut})) + \tanh(\beta \rho_{cut})}
     \\
   & \rho_{cut} \in [0,1] ; \beta > 0
    \end{aligned}
\end{equation}
with $\rho_{cut} = 0.5$ such as not to favor one material. The other parameters $\delta$ and $\beta$ are chosen to preserve the equilibrium between the two parts of the filtering step.

\begin{figure}
\centering
\includegraphics[width=8 cm]{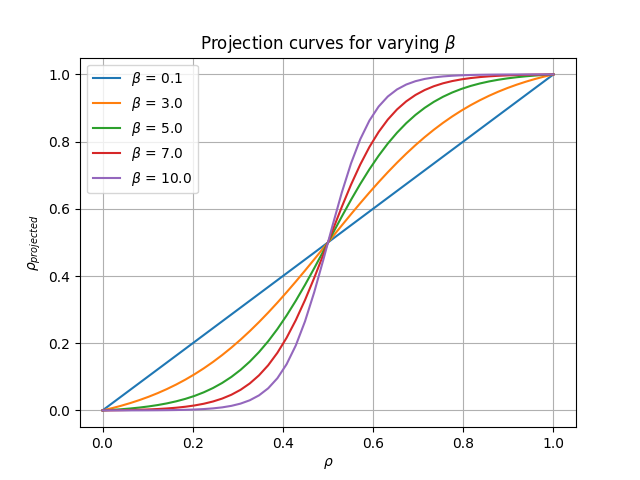}
\caption{Comparison of projection curves for a varying $\beta$ parameter.}
\label{fig:Projection}
\end{figure}

We applied the filtering technique \eqref{eq_filterHelm} and the projection technique \eqref{eq_projection} for the material densities. The same PDE-based filter \eqref{eq_filterHelm} is applied for the densities $\rho_{M_x}$ and $\rho_{M_y}$, to favor uniform magnetization direction in a magnet area while the projection \eqref{eq_projection} is only applied on $|\vec{M}|$ to avoid scaling issues.

\subsection{Direct penalization of intermediate materials}
Some final designs can still present fuzzy boundaries and intermediate material, especially if the optimization starting point is near a local minimum.
To help overcome this issue, we propose to penalize the intermediate materials directly as done in the phase-field topology optimization method \cite{garcke_numerical_2015} and add to the cost function the following term with a weight $\gamma > 0$:
\begin{equation}
    I_{\gamma}(\rho) = \frac{4\gamma}{V_{\Omega_{rot}}}\int_{\Omega_{rot}} \rho(x)(1-\rho(x)) dx.
\end{equation}
The penalization is only applied on iron density $\rho_\nu$ and the magnetization norm $|\vec{M}|$.

\subsection{Post-processing: K-mean heuristic}
In our optimization problem we look for permanent magnetization directions which may change continuously in space. In order to obtain designs which comply with feasibility constraints, we here propose a post-processing step. A K-mean heuristic \cite{macqueen_methods_1967} clustering method is applied. Here, we suggest adapting it to create clusters of elements of similar magnetization direction. \\
Let us define a point $P = (p_x,p_y,p_{\beta_M}) $ by its coordinates $p_x$, $p_y$ in the 2D plane and its magnetization angle $p_{\beta_M}$. We define a set of $k$ points $C_1,\dots C_k$ where $C_j = (c_{x,j}, c_{y,j}, -)$, which we will refer to as centroids, and which are first randomly sampled in the 2D plane.  \\
Let $P_1, \dots, P_{N_{rot}}$ be the centroids of the triangles in the rotor domain $\Omega_{rot}$. At each step of the algorithm, we associate each point $P_i$ with the closest centroid $C_j \in S_k$ using a modified 3D Euclidean distance $d_\alpha$. For two such points $P_i = (p_{x,i}, p_{y,i}, p_{{\beta_M},i})$ and $C_j = (c_{x,j}, c_{y,j}, c_{{\beta_M},j})$, this modified distance function is defined as
\begin{equation}
\begin{aligned}
    d_{\alpha}: & \mathbb R^3 \times \mathbb R^3 \rightarrow \mathbb R \\
    & (P_i,C_j) \mapsto \sqrt{\left(\frac{p_{x,i}-c_{x,j}}{N_x}\right)^2
    +\left(\frac{p_{y,i}-c_{y,j}}{N_y}\right)^2 
    +\alpha \left(\frac{p_{{\beta_M}, i}-c_{{\beta_M},j}}{2\pi}\right)^2}
\end{aligned}
\end{equation}
where $N_x$, $N_y$, $\alpha$ are three weighting constants which can be used to tune the method.

By this procedure, we get $k$ point clusters $S_1, \dots, S_k$ with $S_j = \{ P_i: d_\alpha(P_i, C_j) < d_\alpha(P_i, C_l)$ for all $l \neq j \}$ for $j=1, \dots, k$. Thus, all points in cluster $S_j$ are closer to point $C_j$ in terms of $d_\alpha$ than to any other centroid $C_l$.

The positions of the centroids $C_1, \dots, C_k$ are then updated according to the mean coordinates of the points associated with them, 
\begin{equation}
C_j \leftarrow (\mu_{x,j},\mu_{y,j},\mu_{\beta_{M,j}}) \mbox{ where } \mu_{x,j} = \frac{1}{\#S_j} \sum_{P_i \in S_j} p_{x,i}
\end{equation}
and $\mu_{y,j},\mu_{\beta_{M,j}}$ defined analogously. Here, $\#S_j$ denotes the cardinality of $S_j$.
In the first iteration where $c_{{\beta_M},j}$ is not defined, we use the modified 2D Euclidean distance $d_{\alpha = 0}$.

\section{Application to the magnetostatics problem}
All computations were conducted using the NGSolve \cite{gangl_fully_2020,schoeberl_c11_2014} framework with its python interface.
Three sets of results are presented in this section to validate the proposed strategy. For the different optimization sets, the starting points are given in Table\,\ref{starting}. \\
The first design in Table\,\ref{starting} is used for the Iron-Air optimization.
The results for the material distribution are given in Table\,\ref{IronTable} for different volume fractions $f_{\nu,f}$. These results are coherent with the literature of synchro-reluctant actuators with distributed winding \cite{Vagati1992}, with a phase angle $\varphi = \frac{5\pi}{6}$ and one pair of poles. This phase angle corresponds to the maximum torque of the machine. Results can also be compared with results found in \cite{mellak_synchronous_2018,gangl_multi-objective_2020} which are validated by means of \textit{JMAG} . \\
For the multi-material topology optimization including magnets, the final design obtained at $40\%$ of ferromagnetic material is used as a starting point. Results for magnets distribution are presented in Table\,\ref{MagnetTable}. As expected, magnets are distributed on the air barriers domain. Several constraints on magnets volume were chosen. One should note that the ferromagnetic distribution is modified at the outer radius of the rotor. Without magnets, the reluctance torque is equal to $1.0768$\,[N.m]. Optimized torque and post-processed values after k-mean clustering are given, the number of clusters is fixed to $k=5$. \\
In order to validate the proposed strategy, a new phase angle is chosen. It is equal to $\varphi = \frac{3\pi}{32}$. The objective is to find the optimal materials distribution to maximize the torque at this phase angle. An unbiased starting point is chosen, where homogeneous grey materials are set in the rotor such as $\rho_\nu =0.5,\rho_{M_x} =0.5,\rho_{M_y} =0.5$. The optimal rotor structure is given in Table\,\ref{IronMagnetNeutralTable}. Again, we present optimized torques as well as the torque values after post-processing, where we used $k=5$ clusters for the K-mean clustering. Results on the torque value are comparable for the previously obtained but at higher magnets volume, which is also coherent with literature. On the other hand, the first optimal multi-material model has larger reluctant torque due to optimal current supply in the q-axis.



\begin{table*}
\caption{Starting designs\label{starting}}
\centering
\begin{tabular}{p{0.04\linewidth}p{0.25\linewidth}p{0.25\linewidth}p{0.25\linewidth}}
\hline\noalign{\smallskip}
& & \hfil Starting Design  &\\
\noalign{\smallskip}\hline\noalign{\smallskip}
\hfil \raisebox{7\height}{Rotor} &
 \hfil \includegraphics[width=3.2cm]{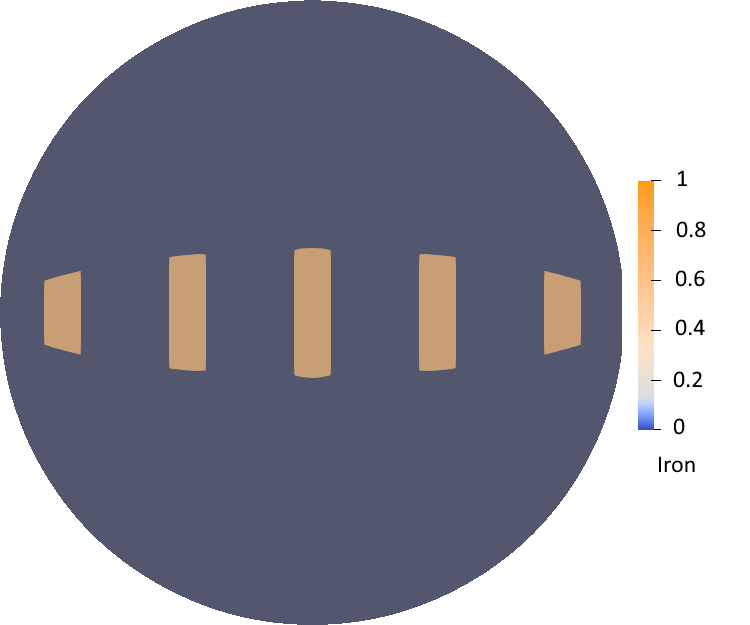} & \hfil \includegraphics[width=3.2cm]{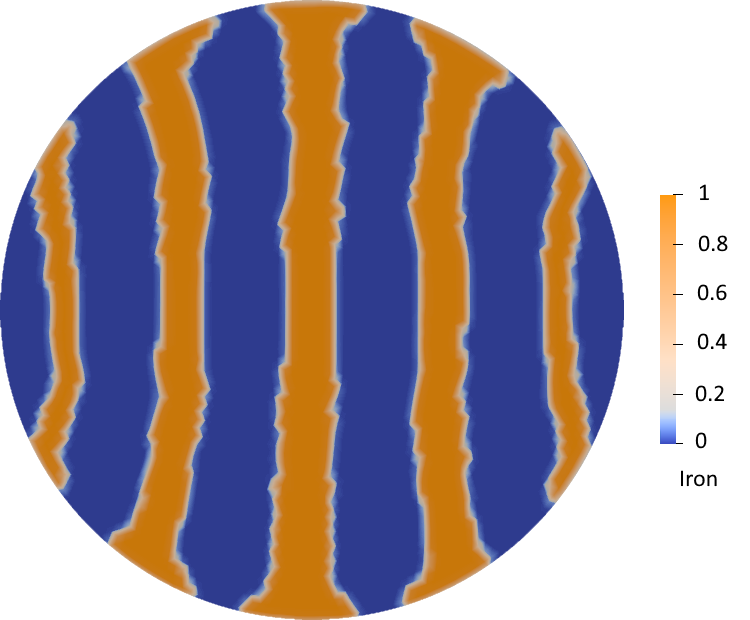} & \hfil  \includegraphics[width=3.2cm]{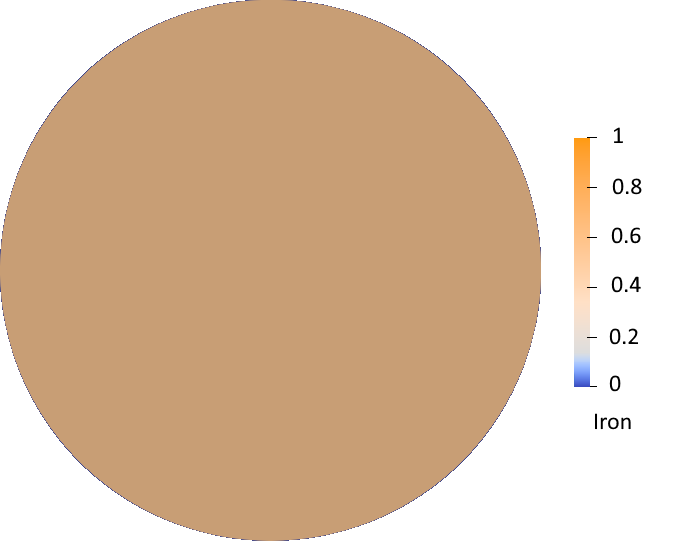}\\
\noalign{\smallskip}\hline

\end{tabular}
\end{table*}

\begin{table*}
\caption{Designs Iron-Air \label{IronTable}}
\centering
\begin{tabular}{p{0.04\linewidth}p{0.25\linewidth}p{0.25\linewidth}p{0.25\linewidth}}
\hline\noalign{\smallskip}
& & \hfil Final filtered design  &\\
\noalign{\smallskip}\hline\noalign{\smallskip}
\hfil \raisebox{7\height}{Rotor} &
 \hfil \includegraphics[width=3cm]{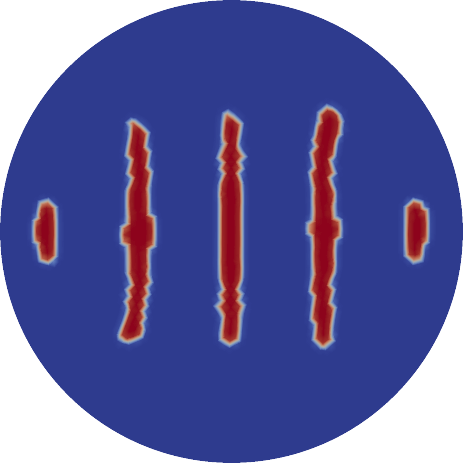} & \hfil \includegraphics[width=3cm]{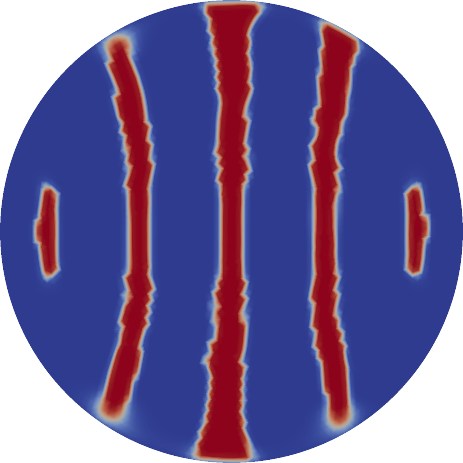} & \hfil \includegraphics[width=3cm]{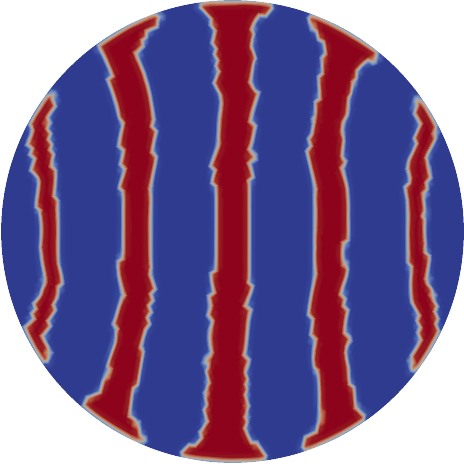}\\
\hfil \raisebox{.5\height}{$f_{\nu,f}$} & \hfil 10 \%  & \hfil 20 \% & \hfil 40 \% \\
\hfil \raisebox{.5\height}{$\bar{T}[N.m]$} & \vfil{\hfil 0.20897}  &  \vfil{\hfil 0.83351} & \vfil{\hfil 1.1129} \\

\noalign{\smallskip}\hline

\end{tabular}
\end{table*}

\begin{table*}
\caption{Designs Magnet-Air-Iron \label{MagnetTable}}
\centering
\begin{tabular}{p{0.04\linewidth}p{0.25\linewidth}p{0.25\linewidth}p{0.25\linewidth}}
\hline\noalign{\smallskip}
& & \hfil Final filtered designs  &\\
\noalign{\smallskip}\hline\noalign{\smallskip}
\hfil \raisebox{7\height}{Rotor} &
 \hfil \includegraphics[width=3cm]{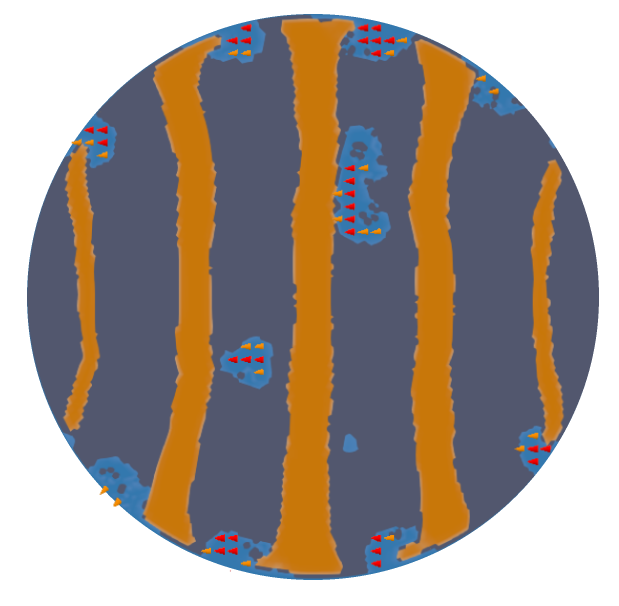} & \hfil \includegraphics[width=3cm]{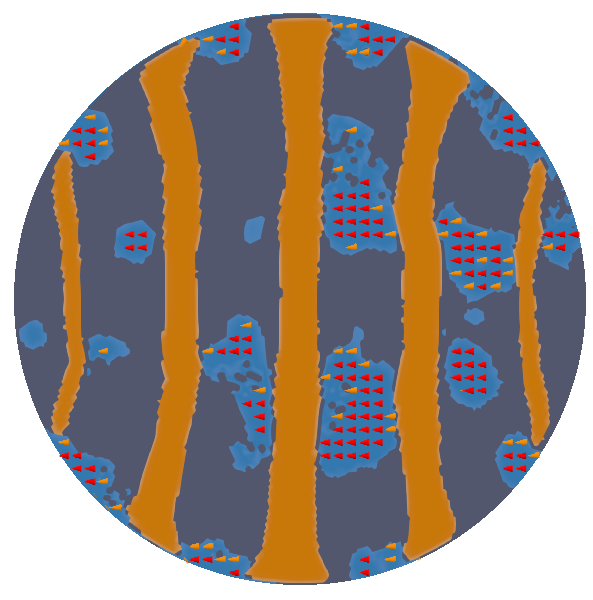} & \hfil \includegraphics[width=3cm]{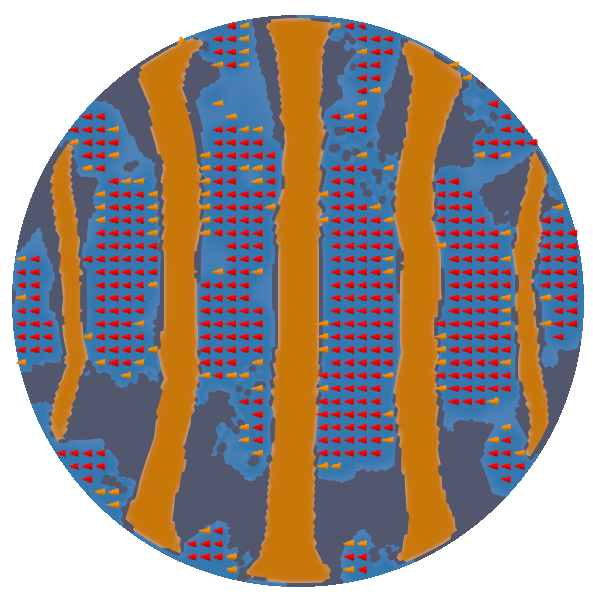}\\
\hfil \raisebox{.5\height}{$f_{\nu,mag}$} & \hfil 7.5 \%  & \hfil 15 \% & \hfil 30 \% \\
\hfil \raisebox{.5\height}{$\bar{T}[N.m]$} & \vfil{\hfil 1.2830}  &  \vfil{\hfil 1.4729} & \vfil{\hfil 1.9000} \\
\hfil \raisebox{.5\height}{$\bar{T}_{Kmeans}[N.m]$} & \vfil{\hfil  1.2646}  &  \vfil{\hfil 1.4469} & \vfil{\hfil 1.8671} \\
\noalign{\smallskip}\hline

\end{tabular}
\end{table*}

\begin{table*}
\caption{Designs Iron-Air-Magnets\label{IronMagnetNeutralTable}}
    \centering
    \begin{tabular}{p{0.04\linewidth}p{0.25\linewidth}p{0.25\linewidth}p{0.25\linewidth}}
    \hline\noalign{\smallskip}
& & \hfil Final Design  &\\
\noalign{\smallskip}\hline\noalign{\smallskip}
\hfil \raisebox{7\height}{Rotor} &
 \hfil \includegraphics[width=3cm]{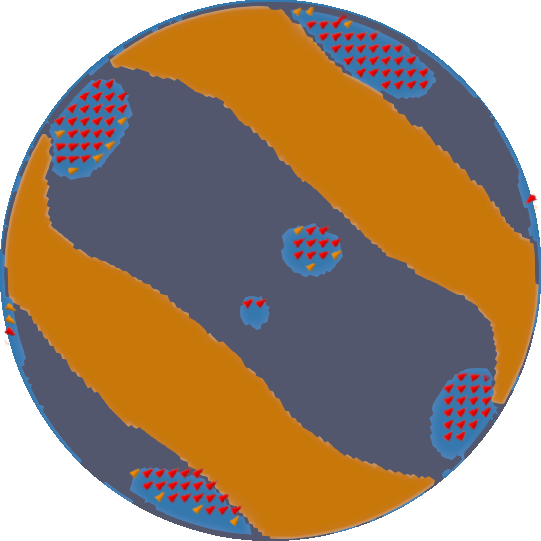} & \hfil \includegraphics[width=3cm]{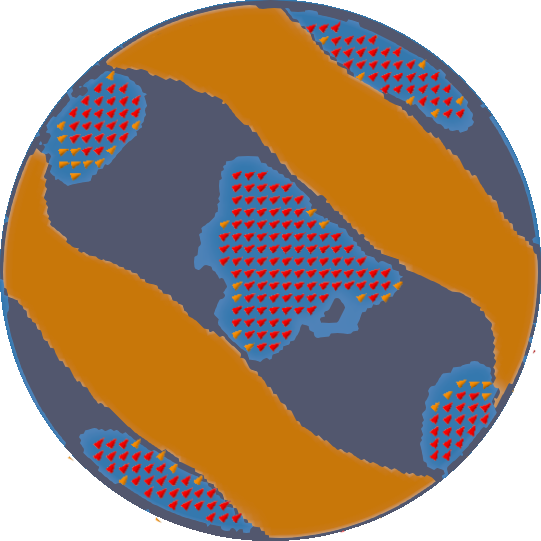}   & \hfil \includegraphics[width=3cm]{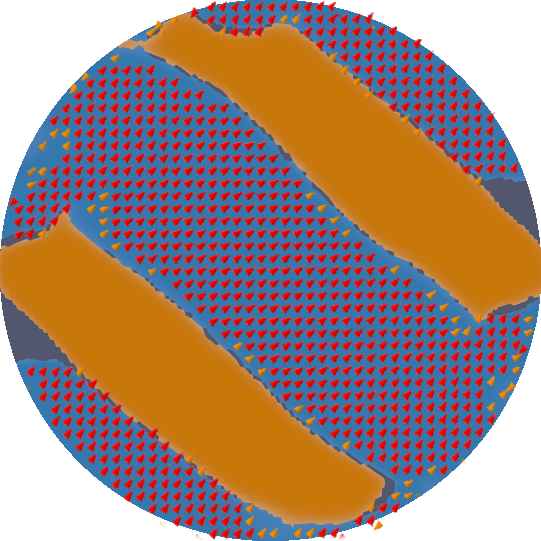}  \\
\raisebox{-1\height}{$f_{v,mag}$} &  \vfil{\hfil 10\%} & \vfil{\hfil 20\%} & \vfil{\hfil Not bounded}  \\
 \raisebox{-1\height}{$\bar{T} \, [N.m]$} & \vfil{\hfil 1.2710} & \vfil{\hfil 1.4513} & \vfil{\hfil 2.0484} \\
 \raisebox{-1\height}{$\bar{T}_{Kmeans} \, [N.m]$} & \vfil{\hfil 1.1097} &\vfil{\hfil 1.3863} &\vfil{\hfil  1.6422} \\
\noalign{\smallskip}\hline
\end{tabular}
\end{table*}

\newpage

\section{Conclusion and outlooks}
In this study we proposed a novel multi-material interpolation method to determine the optimal distribution of air, iron, and magnets for PMSynRM. 
The interpolation takes into account magnetization amplitudes and direction, and a post-processing clustering method is also suggested to homogenize magnets direction for feasibility constraints.  
This study also investigated the use of the four statics positions method to reduce global computation time. 
With further work, an exhaustive comparison with other interpolation schemes will be considered, and the extension to multiobjective optimization under gradient-based methods.
\appendix

\section{Nitsche-mortaring reformulation}
\label{appendix}
Several methods exist to consider the movement in electrical motors, such as the Moving Band (MB) technique described in \cite{Davat1985}, even when using high order elements in the MB, this method remains less accurate than the mortar element method \cite{Antunes2005}. A variant of the mortar method, the Nitsche method, is chosen to take into account the rotation \cite{hollaus_nitsche-type_2010}.
Let us define the adapted operator $K_\theta^{NM}$ to fit the new formulation
\begin{equation}
    \begin{split}
        K_\theta^{NM}(\mathbf{X},u,\eta,v,\mu) = &
        \sum_{i\in \{ rot,stat\}} \left( \int_{\Omega_i}\nu(\rho_\nu^\theta,|\nabla u_i|)\nabla u_i \cdot \nabla v_i \right) \\
        &- \int_{\Omega_{rot}} f_{\nu}(1-\rho_{\nu}^\theta) \frac{M_{max}f_{M}(|\vec{M}|)}{|\vec{M}|}
    R_{\theta} \begin{bmatrix} -M_y^\theta \\ M_x^\theta \end{bmatrix} \cdot \nabla v_{rot} \\
    & + \sum_{i\in \{ rot,stat\}} \left(   
    -\int_{\partial \Omega_i}\nu_0\frac{\partial u_i}{\partial n}(v_i-\mu)
    -\int_{\partial \Omega_i}\nu_0\frac{\partial v_i}{\partial n} (u_i-\eta) \right.) \\
     & \left. \phantom{+\sum_{\Omega\in \{ \Omega_{rot},\Omega_{stat}\}}(}  \; + \frac{\alpha p^2}{h}\int_{\partial \Omega_i}(u_i-\eta)(v_i-\mu) 
    \right),
    \end{split}
\end{equation}
where $(u,\eta,v,\mu) \in V \times W \times V \times W$ and $V,W$ are defined as
\begin{equation}
    \begin{cases}
    V = \{v=(v_{rot}, v_{stat}) \in H^1(\Omega_{rot}) \times H^1(\Omega_{stat}), \, v = 0\text{ on }\partial \Omega \}, \\
    W = \{w \in L^2(\partial \Omega_{rot} \cap \partial \Omega_{stat})\}. \\
    \end{cases}
\end{equation}
Here, $\alpha > 0$ is a stabilization parameter which we chose as $\alpha = 160$, $p=1$ denotes the polynomial degree of the finite element discretization and $h$ the diameter of the largest element of the mesh. 
Moreover, recall the implicit dependence of $\vec{M} = (M_x^\theta, M_y^\theta)^\top$ on the density variables $\rho_{M_x^\theta}, \rho_{M_y^\theta}$ \eqref{eq:def_fsd_tilde}, \eqref{eq:def_fsd_tilde2}.

Hence the state equation \eqref{eq:PDE_density} can be formulated as 
\begin{equation}
\label{eq:stateHNM}
    \mbox{Find } (u_\theta,\eta_\theta) \in V \times W, \; K_\theta^{NM}(\mathbf{X}_\rho, u_\theta,\eta_\theta,v,\mu) = \int_{\Omega_{stat}} j(\theta) \, v dx \quad \mbox{for all } (v,\mu) \in  V \times W.
\end{equation}
In a similar manner, the adjoint equation \eqref{eq:adjoint} can be redefined as
\begin{equation}
\label{eq:adjointHNM}
\mbox{Find } (w_\theta, \mu_\theta) \in V \times W: \frac{\partial K_\theta^{NM}}{\partial (u, \eta)} (\mathbf X,u_\theta, \eta_\theta, w_\theta, \mu_\theta )(\hat u, \hat \eta) = - \frac{1}{4} \frac{\partial T}{\partial u}(u_\theta)(\hat u) 
\end{equation}
for all $(\hat u, \hat \eta) \in V \times W$ for each $\theta \in \{0, \frac{\pi}{12}, \frac{\pi}{6}, \frac{\pi}{4} \}$.

The corresponding Lagrangian reads
\begin{equation}
\tilde {\mathcal L}(\mathbf X, \underline u, \underline \eta, \underline v, \underline \mu) := L (\mathbf X, \underline u) + \sum_{\theta \in \{ 0,\frac{\pi}{12},\frac{\pi}{6},\frac{\pi}{4} \}} K_\theta^{NM}(\mathbf X,  u_\theta, \eta_\theta, v_\theta, \mu_\theta) - \int_{\Omega_{stat}} j(\theta) v_\theta,
\end{equation}
and its gradient is given by
\begin{equation}
\label{eq:diffLagHNM}
    \begin{split}
     \frac{\partial \mathcal{\tilde L}}{\partial \mathbf{X}} (\mathbf X, \underline u, \underline \eta, \underline v, \underline \mu) = &      \left(
     \sum_{\theta \in \{0,\frac{\pi}{12},\frac{\pi}{6},\frac{\pi}{4}\}} \frac{\partial K^{NM}_\theta}{\partial \mathbf{X}}(\mathbf X , u_\theta,\eta_\theta,v_\theta,\mu_\theta)  -\frac{1}{4} \frac{\partial T_\theta}{\partial \mathbf{X}}\right) \\
     & + 
     \frac{\partial \psi(h_{v,f}(\mathbf{X}),\gamma_f,\mu)}{\partial \mathbf{X}} +
     \frac{\partial \psi(h_{v,mag}(\mathbf{X}),\gamma_{mag},\mu)}{\partial \mathbf{X}}.
    \end{split}
\end{equation}

\vspace{6pt} 
\bibliographystyle{abbrv}
\bibliography{bib.bib}

\end{document}